\documentclass[12pt,a4paper]{article}
\usepackage{amssymb}

\usepackage{amssymb}
\usepackage{hyperref}

\usepackage{graphicx,amssymb,amsfonts,epsfig,a4,amsmath,amsthm,xypic}
\usepackage[latin1]{inputenc}

\usepackage{tikz-cd}

\newtheorem{Thm}{Theorem}[section]
\newtheorem{Prop}[Thm]{Proposition}
\newtheorem{Lem}[Thm]{Lemma}
\newtheorem{Cor}[Thm]{Corollary}

\newtheorem{Rem}[Thm]{Remark}

\newtheorem{quest}{Question}

\newcommand{\Z}{\mathbb{Z}}

\newcommand{\R}{\mathbb{R}}

\newcommand{\Q}{\mathbb{Q}}

\newcommand{\bpr}{\noindent \textbf{Proof}: }

\newcommand{\epr}{~$\blacksquare$}

\newcommand{\ra}{\rightarrow}

\title{On the transportation cost norm on finite metric graphs}
\author{Georges SKANDALIS and Alain VALETTE}
\date{January 23, 2026}

\begin{document}

\maketitle

\begin{abstract} For a finite metric graph $X=(V,E,\ell)$, where $V$ is endowed with the shortest path metric, we consider the transportation cost problem associated with the distance $d$ on $V$. Namely, for $f$ a function with total sum 0 on $V$, write $f=\sum_{a,b\in V}P(a,b)(\delta_a-\delta_b)$ where the transportation plan $P$ satisfies $P(a,b)\geq 0$ for $(a,b)\in V\times V$. The cost of $P$ is $W(P):=\sum_{a,b\in V}P(a,b)d(a,b)$ and the transportation norm of $f$ is $\|f\|_{TC}=\min_P W(P)$ where $P$ runs over all transportation plans for $f$.

In this semi-survey paper, we give short proofs for the following statements:
\begin{itemize}
\item There always exists an optimal transportation plan supported in $V_+\times V_-$ where $V_+=\{x\in V: f(x)>0\}$ and $V_-=\{x\in V: f(x)<0\}$. If $X$ is a metric tree, we may moreover assume that this plan involves at most $|Supp(f)|-1$ transports. 
\item There always exists an optimal transportation plan supported in the set of edges of $X$.
\item Better, there always exists an optimal transportation plan supported in some spanning tree of $X$. 
\end{itemize}
We use this to reprove known formulae for the transportation norm when $X$ is either a tree or a cycle.
\end{abstract}

\section{Introduction}

Let $(X,d)$ be a finite metric space. The {\it transportation cost space} of $X$ (also known as {\it Arens-Eells space, Kantorovich-Rubinstein space, or Lipschitz free space}\footnote{See Definition 1.1 in \cite{GodKal}, or Chapter 3 in \cite{Weaver}; for historical background see the remarks in section 10.1 of \cite{Ostro} and section 1.6 in \cite{OstOst1}.}) is the space $\mathcal{F}(X)$ of real-valued functions $f$ on $X$ with total mass 0, endowed with a norm defined as follows. For such an $f$, set 
$$\Pi_f=\{P:X\times X\rightarrow\R_{\geq 0}: \;f=\sum_{a,b\in X} P(a,b)(\delta_a -\delta_b)\},$$
where $\delta_x$ denotes the Dirac mass at $x$. Equivalenty $P\in\Pi_f$ if and only if 
$$f(x)=\sum_{a\in X}(P(x,a)-P(a,x))$$
for all $x\in X$. It is clear that $\Pi_f$ is non-empty, and does not depend on $d$. Put also $X_+=\{a\in V;\ f(a)>0\}$ and $X_-=\{a\in V;\ f(a)<0\}$ and let also $$\widetilde \Pi_f=\{P\in \Pi_f;\ \mathrm{supp}(P)\subset X_+\times X_-\}.$$  It is clear that $\widetilde \Pi_f$ is also non-empty. The elements of $\Pi_f$ are \emph{transportation plans} for $f$ and those of $\widetilde \Pi_f$ are \emph{simultaneous transportation plans.}

We take $d$ into account by defining the {\it cost} of $P\in\Pi_f$ as:
$$W_X(P)=\sum_{a,b\in X} P(a,b)d(a,b).$$
Now define the {\it transportation cost norm} of $f$ as
$$\|f\|_{TC(X)}=\min\{W(P):P\in\Pi_f\}.$$
and the {\it simultaneous transportation cost norm} of $f$ as
$$\|f\|_{STC(X)}=\min\{W(P):P\in\widetilde\Pi_f\}.$$

Let $P_0\in \Pi_f$. The set $\{P\in \Pi_f; \ W(P)\le W(P_0)\}$ is compact, and therefore the infimum is a minimum. When the metric space $X$ is clear from context we drop the subscript, both in $W(P)$ and in $\|f\|_{TC}$.


We recall the optimal transportation interpretation of $\|f\|_{TC}$: suppose that a mass $f(x)$ is available at points $x\in X$ where $f(x)>0$, and a mass $-f(y)$ is in demand at points $y\in X$ with $f(y)<0$. An element $P\in \Pi_f$ can be viewed as a {\it transportation plan}, i.e. a way of transporting the mass from points where it is in excess to the points where it is in demand, $W(P)$ represents the physical work associated with this transportation plan, and $\|f\|_{TC}$ represents the most economical transportation associated with $f$. 

The following result is due to Bergman (Proposition 24 in \cite{Berg08}) with a short proof in Proposition 3.16 of \cite{Weaver}. We will give a new proof in Section 2.
\begin{Thm}\label{ThmSTC}
Given any finite metric space $(X,d)$ and any $f\in TC(X)$, we have $\|f\|_{TC}=\|f\|_{STC}$.
\end{Thm}



Next we deal with $X=(V,E,\ell)$ a finite, connected, metric graph. This means that $(V,E)$ is a simple connected graph with no loop, $\ell:E\rightarrow \R_{>0}$ is a function giving a positive length to every edge, and $V$ becomes a finite metric space with the shortest path metric. We define two more transportation cost norms on $\mathcal{F}(V)$ associated with the graph structure.

\begin{itemize}
\item Let $\sim$ denote the adjacency relation on $V$, so that $x\sim y$ if and only if $\{x,y\}\in E$. Let $\Pi_{f,E}$ denote the set of transportation plans supported in the graph of $\sim$ in $V\times V$. This corresponds to transports supported on the edges of $X$. We denote by $\|f\|_{TCE}$ the corresponding norm.
\item Recall that a {\it spanning tree} in $X$ is a subtree of $X$ passing through all vertices. We denote by $\Pi_{f,T}$ the set of transportation plans $P\in\Pi_f$ for which there exists a spanning tree $T$ such that $P$ is supported on the edges of $T$. We denote by $\|f\|_{TCT}$ the corresponding norm.
\end{itemize}

Assume that $X$ is endowed with an orientation, i.e. every edge $e\in E$ has an origin $e_-$ end an extremity $e_+$. This allows to define the boundary operator $\partial:\R^E\rightarrow\R^V:\delta_e\mapsto\delta_{e_+}-\delta_{e_-}$; equivalently
\begin{equation}\label{bord}
(\partial\phi)(v)=\sum_{e:e_+=v}\phi(e)-\sum_{e:e_-=v}\phi(e)
\end{equation}
($v\in V, \phi\in\R^E$). It is classical that the image of $\partial$ is exactly the space $\mathcal{F}(V)$, so that the rank of $\partial$ is $|V|-1$; and that the kernel of $\partial$ is the space $Z(X)$ of cycles of $X$, so that $\dim(Z(X))=|E|-|V|+1$ (see Proposition 4.3 and Theorem 4.5 in \cite{Biggs}). Then:

\begin{Thm}\label{ostost}
\begin{enumerate}
\item For $f\in\mathcal{F}(V)$ we have:
$$\| f\|_{TC}=\|f\|_{TCE}=\|f\|_{TCT}.$$
\item The operator $\partial:\R^E\rightarrow\R^V$ induces an isometric isomorphism from the quotient space $\ell^1(E,\ell)/Z(X)$ onto $TC(V)$. 
\hfill$\blacksquare$
\end{enumerate}
\end{Thm}

The first equality in the first part is proved after Remark 1.3 in \cite{OstOst3}; the 2nd equality might be new. The latter corresponds to the intuitive idea that, if a transportation plan involves loops in the graph, it is possible to avoid the loops without increasing the cost.

For the second part, see Proposition 10.10 in \cite{Ostro} for usual graphs, and Theorem 3.1 in \cite{OstOst2} or Proposition 1.5 in \cite{OstOst3} for the general case; note that that the proofs in \cite{Ostro} and \cite{OstOst2} are indirect, establishing that the transpose operator $\partial^t:(\mathcal{F}(V))^*\rightarrow (\ell^1(E,\ell)/Z(X))^*$, is an isometric isomorphism. A short direct proof of Theorem \ref{ostost} will be given in Section 3.

Section 4 will be devoted to 3 corollaries of Theorem \ref{ostost}.

\begin{itemize}

\item Due to the presence of a minimum either in the original definition of $\|f\|_{TCE}$, or as a quotient norm in Theorem \ref{ostost}, for a general metric graph we don't expect a closed formula for $\|f\|_{TCE}$. A notable exception to this philosophy is the case of metric {\it trees}: then an explicit formula for $\|f\|_{TC}$ is available. We refer to  \cite{MPV} for the interesting history and two proofs of the formula; see also section 2.1 in \cite{Men}.

\begin{Cor}Let $X=(V,E,\ell)$ be a finite metric tree. Fix an arbitrary orientation of the edges. For $e\in E$, we define $S_e$ as the vertex set of the half-tree containing $e_+$, i.e. $S_e=\{x\in V; \ d(x,e_+)\le d(x,e_-)\}$.   For $f\in\mathcal{F}(V)$ we have:
\begin{equation}\label{magic} 
\|f\|_{TCE}=\sum_{e\in E}|f(S_e)|\ell(e).
\end{equation}
wehere $f(S_e)=\sum_{x\in S_e}f(x)$.
\end{Cor}

Our new, simple proof of (\ref{magic}) follows by combining Theorem \ref{ostost} with injectivity of $\partial$ for trees - see Proposition \ref{cortrees} below.

\item In Corollary \ref{cycle} we provide a new proof of a formula for the transportation cost norm on the $N$-cycle $C_N$, originally appearing as formula (4.2.) in \cite{CaMo}.

\item Recall that a {\it bridge} in a graph is an edge whose deletion produces a graph with more connected components\footnote{Equivalently, an edge is a bridge if it is not contained in any cycle in the graph.}. Deleting all bridges in a graph results in a {\it bridgeless} graph. We show in Corollary \ref{robust} that, for a general finite metric graph $X$, solving the TC problem is equivalent to solving it on the bridgeless part $X_b$ of $X$.

\end{itemize}

Finally Section 5 discusses some algorithmic issues in the case of metric trees. For example, we provide an algorithmic way to find $P\in\tilde{\Pi}_f$ such that $W(P)=\|f\|_{STC}$ and $|Supp(P)|\le |Supp(f)|-1$.

\section{A new proof of Theorem \ref{ThmSTC}}

The proof will follow from the following Lemma.

\begin{Lem}\label{strict}
If all triangle inequalities in $X$ are strict \emph{i.e.}, if for all $x,y,z$ in $X$ with $x\ne y\ne z$, we have $d(x,z)<d(x,y)+d(y,z)$ and $W(P)=\|f\|_{TC}$, then $\mathrm{supp}(P)\subset X_+\times X_-$.
\end{Lem}
\bpr
If there exist $(x,y)\in \mathrm{supp}(P)$ with $y\not\in X_-$, then $0\leq f(y)=\sum_{a\in X}(P(y,a)-P(a,y))$ and, since $P(x,y)>0$, there exists $z\in X$ with $P(y,z)> 0$. Let $m=\min(P(x,y),P(y,z))$ and define $Q$ by setting $Q(x,z)=m+P(x,z)$, $Q(x,y)=P(x,y)-m$, $Q(y,z)=P(y,z)-m$, and $Q(u,v)=P(u,v)$ for all other points. It is readily checked that $Q\in \Pi_f$ and $W(P)-W(Q)=a(d(x,y)+d(y,z)-d(x,z))$, contradicting $W(P)=\|f\|_{TC}$.

We do a similar reasoning if there exist $(x,y)\in \mathrm{supp}(P)$ with with $x\notin X_+$.\epr

\bigskip
{\bf Proof of Theorem \ref{ThmSTC}:} For $n\geq 1$, define a distance $d_n$ on $X$ by $d_n(x,x)=0$ and $d_n(x,y)=\frac{1}{n}+d(x,y)$ for $x\ne y$. Write $\|f\|_{TC_n}$ for the transportation cost norm associated with $d_n$. Let $P_n\in\Pi_f$ minimize the cost for $d_n$. By lemma \ref{strict}, $P_n$ is supported in $X_+\times X_-$. By compactness, we may assume that $(P_n)_{n>0}$ converges pointwise to $P$, also supported in $X_+\times X_-$. As $(d_n)_{n>0}$ converges pointwise to $d$, we also have $W(P)=\lim_{n\ra\infty}W_n(P_n)=\lim_{n\ra\infty}\|f\|_{TC_n}$.

It remains to see that $W(P)=\|f\|_{TC}$. To see this, let $Q\in\Pi_f$ be such that $W(Q)=\|f\|_{TC}$. Then
$$\|f\|_{TC_n}\leq W_n(Q)\leq W(Q)+\frac{1}{n}\|Q\|_1=\|f\|_{TC}+\frac{1}{n}\|Q\|_1$$
hence, passing to the limit: $W(P)\leq\|f\|_{TC}$. The other inequality holds by definition. $\blacksquare$



\section{Towards Theorem \ref{ostost}}

\subsection{The basic construction}

For every 2-element subset $\{a,b\}\subset V$, we select (arbitrarily if non-unique) a shortest length edge path $[a,b]$ joining $a$ and $b$. Let $\varepsilon_{ab}\in \R^E$ be the {\it characteristic function of $[a,b]$}, defined as:
$$\varepsilon_{ab}(e)=\left\{\begin{array}{cccc}0 & \mbox{if}  & e\notin[a,b]; & \\1 & \mbox{if} & e\in [a,b]\;\mbox{and} &e\;\mbox{points from}\;a\;\mbox{to}\;b; \\-1 & \mbox{if} & e\in [a,b] \;\mbox{and} & e \;\mbox{points from}\; b\;\mbox{to}\;a.\end{array}\right.$$

To every $P:V\times V\rightarrow\R_{\geq 0}$, we associate $\phi_P\in\R^E$ defined by:
$$\phi_P=\sum_{a,b\in V}P(a,b)\varepsilon_{ab}.$$
Motivation for introducing $\phi_P$ comes from the following simple but crucial

\medspace

{\bf Basic observation:} For a fixed $f\in \mathcal{F}(V)$ and $P\geq 0$, we have $P\in\Pi_f$ if and only if $\partial\phi_P=-f$.

\medspace

Indeed $\partial\phi_P=\sum_{a,b\in V}P(a,b)(\delta_b-\delta_a)$, so it follows from definition of $\Pi_f$.

\medskip
\begin{Lem}\label{easyineq} For every $P:V\times V\rightarrow\R_{\geq 0}$ we have:
$$\|\phi_P\|_1\leq W(P),$$
where $\|\cdot\|_1$ denotes the norm in $\ell^1(E,\ell)$.
\end{Lem}

\bpr Observe that $d(a,b)=\sum_{e\in E}\ell(e)|\varepsilon_{ab}(e)|$ for every $a,b\in V$, hence 
$$W(P)=\sum_{a,b\in V}\sum_{e\in E} P(a,b)\ell(e)|\varepsilon_{ab}(e)|=\sum_{e\in E}\ell(e)\sum_{a,b\in V} P(a,b)|\varepsilon_{ab}(e)|$$
$$\geq \sum_{e\in E}\ell(e)\,\Big|\sum_{a,b\in V}P(a,b)\varepsilon_{ab}(e)\Big|=\sum_{e\in E}\ell(e)|\phi_P(e)|=\|\phi_P\|_1.$$
\epr

\begin{Lem}\label{onto} For any $f\in TC(V)$ and $\phi\in\partial^{-1}(-f)$, there exists $P\in\Pi_{f,E}$ such that $\phi_P=\phi$ and $W(P)=\|\phi\|_1$.
\end{Lem}

\bpr 
We define $\phi$ as follows:
\begin{itemize}
\item $P(a,b)=\phi(e)$ if $a=e_-,b=e_+$ and $\phi(e)>0$;
\item $P(a,b)=-\phi(e)$ if $a=e_+,b=e_-$ and $\phi(e)<0$;
\item $P(a,b)=0$ in all other cases.
\end{itemize}
Observe that $P(e_-,e_+)-P(e_+,e_-) =\phi(e)$ for every $e\in E$. Then, as $P$ is supported on edges, we have
$$\phi_P=\sum_{a,b\in V}P(a,b)\varepsilon_{ab}=\sum_{e\in E}(P(e_-,e_+)-P(e_+,e_-))\delta_e=\sum_{e\in E}\phi(e)\delta_e=\phi.$$
Then $\partial\phi_P=-f$ so $P \in\Pi_{f,E}$ thanks to the basic observation. Finally
$W(P)=\sum_{a,b\in V}P(a,b)d(a,b)=\sum_{e\in E}|\phi(e)|\ell(e)=\|\phi\|_1.$
\epr

\bigskip

\subsection{Proof of Theorem \ref{ostost}, but the equality $\|f\|_{TC}=\|f\|_{TCT}$}

As $\partial$ realizes a linear isomorphism $\ell^1(E,\ell)/Z(X)\rightarrow\mathcal{F}(V)$ we may endow $\mathcal{F}(V)$ with the quotient norm $\|\cdot\|_{\ell^1(E,\ell)/Z(X)}$. For $f\in\mathcal{F}(V)$ we compute:
$$\|f\|_{\ell^1(E,\ell)/Z(X)}=\inf_{\phi:\partial\phi=-f}\|\phi\|_1\leq\inf_{P\in\Pi_f}\|\phi_P\|_1\leq\inf_{P\in\Pi_f}W(P)=\|f\|_{TC}\leq\|f\|_{TCE}$$
where the first inequality follows from $\partial\phi_P=-f$ for $P\in\Pi_f$, and the second inequality follows from lemma \ref{easyineq}.

Now, take $\phi$ such that $\partial\phi=-f$ and $\|f\|_{\ell^1(E,\ell)/Z(X)}=\|\phi\|_1$. By lemma \ref{onto}, we find $P\in\Pi_{f,E}$ such that $\phi=\phi_P$ = $\|\phi\|_1=W(P)$. This shows that $\|f\|_{TCE}\leq \|f\|_{\ell^1(E,\ell)/Z(X)}$. At this point we proved $\|f\|_{\ell^1(E,\ell)/Z(X)}=\|f\|_{TC}=\|f\|_{TCE}$, so that Theorem \ref{ostost} is proved except for $\|f\|_{TC}=\|f\|_{TCT}$.\hfill$\blacksquare$


 

\subsection{Spanning trees and the proof of $\|f\|_{TC}=\|f\|_{TCT}$}

\begin{Lem}\label{spanning} For every $f\in\mathcal{F}(V)$, there exists $\phi\in\partial^{-1}(f)$ realizing the minimum of $\{\|\psi\|_1:\partial\psi=f\}$ and such that, for every cycle $C$ in $X$, the function $\phi$ vanishes on at least one edge of $C$; equivalently, $\phi$ is supported in some spanning tree $T$ of $X$.
\end{Lem}

\bpr Take any $\phi_0$ realizing the minimum of the $\ell^1$-norm on $\partial^{-1}(-f)$. Let $C$ be any cycle in $X$; recall the definition of the characteristic function $\xi_C$ (taken from p.23 of \cite{Biggs}): choose one cyclic orientation for $C$ and define, for $e\in E$:
$$\xi_C(e)=\left\{\begin{array}{lll}0 & \mbox{if}  & e\notin C; \\1 & \mbox{if} & e \in C \;\mbox{and $e$ coherent with the orientation of $C$} ;\\-1 &\mbox{ if} & e \in C \;\mbox{and $e$ not coherent with the orientation of $C$}. \end{array}\right.$$
Define then $\phi_t=\phi_0+t\xi_C$, for $t\in \R$. Observe that $\phi_t\in\partial^{-1}(-f)$ as $\partial\xi_C=0$. The function $t\mapsto\|\phi_t\|_1=\sum_{e\in E}\ell(e)|(\phi_0+t\xi_C)(e)|$ is convex and piecewise affine; it is affine at every $t$ such that $\phi_t$ does not vanish on $C$. So the minimum of $\|\phi_t\|_1$ is attained for a $t_0$ such that $\phi_{t_0}$ vanishes at some edge in $C$. Note that $\|\phi_{t_0}\|_1=\|\phi_0\|_1$ by minimality of $\|\phi_0\|_1$.

Choosing now $\phi\in\partial^{-1}(-f)$ realizing the minimum of the $\ell^1$-norm and with minimal support, the above argument shows that the support of $\phi$ cannot contain any cycle, i.e. it is a forest, which is of course contained in some spanning tree. 
\epr

\bigskip
{\bf End of proof of Theorem \ref{ostost}}: For $T$ a spanning tree in $X$, we denote by $d_T$ (resp. $\partial_T$) the distance (resp. boundary operator) of $T$. Fix $f\in\mathcal{F}(V)$. For $a,b\in V$, we have $d(a,b)\leq d_T(a,b)$, hence $W_X(P)\leq W_T(P)$ for every $P\in \Pi_f$, and $\|f\|_{TC(X)}\leq \|f\|_{TC(T)}$. To show that the equality is realized by some spanning tree, by the part already proved of Theorem \ref{ostost} we find $\phi\in\ell^1(E,\ell)$ such that $\partial\phi=-f$ and $\|f\|_{TC(X)}=\|\phi\|_1$. By lemma \ref{spanning} we may assume that $\phi$ is supported on some spanning tree $T$ in $X$, so that the $\ell^1$-norms of $\phi$ on $E$ and $E(T)$ are equal. Since $\partial_T\phi=-f$, we have $\|f\|_{TC(X)}=\|\phi\|_1=\|f\|_{TC(T)}$.
\hfill$\blacksquare$


\begin{Rem}\label{rem3.4}
Let us come to the general case, where $(X,d)$ is a finite metric space. Let $f\in \mathcal{F}(V)$. By theorem \ref{ThmSTC}, we already know that $\|f\|_{TC}$ doesn't change if we replace $(X,d)$ by $(Supp(f),d)$. Consider the complete  graph with vertex set $Supp(f)$ whose path length is given by $d$. Theorem  \ref{ostost} then implies that there is a spanning tree realizing the optimal transportation, in particular at most $|Supp(f)|-1$  transports are needed.

Note that the number $|Supp(f)|-1$  is generically optimal since, writing $f=\partial (P)$, we see that the $\Q$-rank of the $\Q$-linear span of $\{f(u);\ u\in V\}$ is at most $|Supp (P)|$.
\end{Rem}

\begin{Rem}
In the proof of Lemma \ref{spanning}, we see that the derivative of $t\mapsto \|\varphi_t\|_1$ at any regular point is $\sum _e\pm\xi_C(e) \ell(e)$. In a generic case (in particular whenever all the $\ell(e)$'s are independent over the rational numbers), this derivative cannot be $0$ and it follows that $t\mapsto \|\varphi_t\|_1$ has a unique minimum at $t=0$, so that $\phi_0$ is automatically carried by a spanning tree. We can then give another proof of the equality $\|f\|_{TC}=\|f\|_{TCT}$ in the spirit of our proof of Theorem \ref{ThmSTC}: approximate $\ell$ by a sequence $\ell_n$ such that for all $n$ the family $(\ell_n(e))_{e\in E}$ is free over $\Q$.
\end{Rem}

\section{Corollaries of Theorem \ref{ostost}}

\subsection{The case of trees}

We improve Formula (\ref{magic}) by finding an explicit function $\phi_f$ with $\partial(\phi_f)=-f$ and an explicit $P\in\Pi_f$ with $W(P)=\|f\|_{TC}$.


\begin{Prop}\label{cortrees} Let $X=(V,E,\ell)$ be a finite metric tree, and $f\in \mathcal{F}(V)$.
\begin{enumerate}
\item Set $\phi_f(e):=-f(S_e)=-\sum_{x\in S_e}f(x)$; then $\phi_f$ is the unique solution of $\partial\phi=-f$.
\item We have $$\|f\|_{TCE}=\|\phi_f\|_1=\sum_{e\in E}|f(S_e)|\ell(e)$$
(which is exactly Formula (\ref{magic})).
\item Define $P:V\times V\rightarrow\R^+$ by $P(a,b)=|f(S_e)|$ if either $e_-=a,e_+=b$ and $f(S_e)<0$; or $e_-=b,e_+=a$ and $f(S_e)>0$; and $P(a,b)=0$ in all other cases. Then $P\in\Pi_f$ and $W(P)=\|f\|_{TC}$.
\end{enumerate}\end{Prop}

\bpr \begin{enumerate}
\item It is readily checked that $\partial(\phi_f)=-f$, and moreover $\partial $ is injective as $X$ is a tree.
\item The formula immediately follows from Theorem \ref{ostost}.
\item If we apply the proof of lemma \ref{onto} to $\phi_f$, we find that $\phi_f=\phi_P$, so that $P\in\Pi_f$ and $W(P)=\|\phi_f\|_1$.\epr

\end{enumerate}

\subsection{The case of cycles}



Here we provide a new proof of a formula for the transportation cost norm on the $N$-cycle $C_N$, originally appearing as formula (4.2.) in \cite{CaMo} (see also section 2.2 in \cite{Men} for another proof).

\begin{Cor}\label{cycle} Let $C_N$ be the $N$-cycle with vertex set $\{v_i;\ i\in\Z/N\Z\}$ and edge set $E=\{e_i;\ i\in\Z/N\Z\}$, with the edge $e_i$ joining $v_i$ and $v_{i-1}$. Let $\ell:E\to \R_+^*$ be a function. For $f\in\mathcal{F}(V)$, and $k\in\{0,\ldots,N-1\}$ set $\alpha_k=\sum_{i=0}^k f(i)$. Then 
$$\|f\|_{TC}=\min_{0\leq k\leq N-1}\sum_{i=0}^{N-1}|\alpha_i-\alpha_k|\ell(e_i) $$
This minimum is $\alpha_k$ such that  $$\sum _{i\in \{0,\ldots,N-1\};\ \alpha_i<\alpha_k}\ell(e_i)\le \dfrac12\sum _{i=0}^{N-1}\ell (e_i)\le \sum _{i\in \{0,\ldots,N-1\};\ \alpha_i\le\alpha_k}\ell(e_i).$$ 
In particular, if $\ell(e)=1$ for every vertex, we find $$\|f\|_{TC}=\sum_{i=0}^{N-1}|\alpha_i-\mbox{Me}|, $$
where $\mbox{Me}$ denotes the median of the $\alpha_k$'s.
\end{Cor}

\bpr Choose the orientation of $E$ by setting $(e_i)_+=v_i$ and $(e_i)_-=v_{i-1}$ and put $\phi(e_i)=\alpha_i$, so that $\partial \phi=-f$. The cycle space of $C_N$ is 1-dimensional and generated by the constant function 1 on $E$. By Theorem \ref{ostost}, we know that $\|f\|_{TC}$ is the minimum of $q(t)=\|\phi-t\|_1=\sum_{i}|\alpha_i-t|\ell(e_i)$. For $t\not \in \{\alpha_k;\ k\in \Z/N\Z\}$, we have $$q'(t)=\sum_{i:\ \alpha_i<t}\ell(e_i)-\sum_{i:\ \alpha_i >t}\ell(e_i).$$


By convexity, the minimum of $q(t)$ is attained at $t=\alpha_k$ where $k$ is characterized by $q'(t)\le 0$ for $t<\alpha_k$ and $q'(t)\ge 0$ for $t\ge \alpha_k$. The conclusion follows. \epr

\begin{Rem}
Assume $X$ has just one cycle. Write $V=V'\sqcup \{v_i,\ i\in \Z/N\Z\}$ and $E=E'\sqcup \{e_i,\ i\in \Z/N\Z\}$.
Note that $E'$ is the bridgeless part of the graph.

 Let $f\in \mathcal{F}(V)$. Choose an orientation of $E$ such that $(e_i)_+=v_i$ and $(e_i)_-=v_{i-1}$. Take $\phi\in \ell^1(E,\ell)$ be such that $\partial \phi=-f$ (\emph{e.g.} by removing $e_0$). All the solutions are then $\phi-t\psi$ where $\psi(e_i)=1$ and $\psi(e)=0$ if $e\in E'$. The minimum is obtained with $t=\phi(e_k)$ where $\phi(e_k)$ is such that $$\sum _{i\in \{0,\ldots,N-1\};\ \phi(e_i)<\phi(e_k)}\ell(e_i)\le \dfrac12\sum _{i=0}^{N-1}\ell (e_i)\le \sum _{i\in \{0,\ldots,N-1\};\ \phi(e_i)\le \phi(e_k)}\ell(e_i).$$ 
\end{Rem}

\subsection{An application to arbitrary graphs} 


Fix an orientation of $E$; recall that, for $e\in E$, we defined $S_e=\{x\in V:d(x,e_+)<d(x,e_-)\}$.

\begin{Cor}\label{robust} Solving the TC problem on $X$ is equivalent to solving it on the bridgeless part $X_b$.
\end{Cor}

\bpr Fix $f\in\mathcal{F}(V)$ and $\phi\in\ell^1(E)$ such that $\partial\phi=-f$. Let $B$ denote the set of bridges of $X$.

{\bf Claim:} For $e^0\in B$ we have $\phi(e^0)=f(S_{e^0})$.

Indeed, since $e^0$ disconnects $X$, we may write $E=\{e^0\}\cup E_+\cup E_-$ where $E_+$ is the set of edges whose endpoints are both in $S_{e^0}$ and $E_-$  is the set of edges whose endpoints are both in $V\setminus S_{e^0}$.

Put $\psi=\mathbf{1}_{E_+\cup \{e^0\}}\phi$. We see that $\partial \psi $ coincides with $\partial \phi=-f$ on $S_{e^0}$ and vanishes outside $S_{e^0}\cup e^0_-$. We find $0=\sum_{v\in V}\partial \psi(v)=\sum_{v\in S_{e^0}}\partial \phi(v)+\partial \psi(e^0_-)$. But $\partial \psi(e^0_-)=\psi(e^0)=\phi(e^0)$, which proves the Claim.


\medskip
As a consequence of the Claim, if $\partial\phi=-f$, then 
$$\|\phi\|_1=\sum_{e\in B} |f(S_e)|\ell(e) + \sum_{e\notin B}|\phi(e)|\ell(e),$$
so by Theorem \ref{ostost} we have to minimize the quantity $\sum_{e\notin B}|\phi(e)|\ell(e),$.

We now show that this amounts to solving a TC problem on the bridgeless part $X_b$. For this, define $\phi_0\in\ell^1(E,\ell)$ by 
$$\phi_0=\left\{\begin{array}{ccc}f(S_e) & if & e\in E_d \\0 & if & e\notin E_d\end{array}\right.$$
and set $f_0:=f-\partial\phi_0$. Clearly the equation $\partial\phi=-f_0$ is equivalent to $\partial(\phi-\phi_0)=-f$. But if $\partial\phi=-f_0$ then, by the Claim, for $e\in B$ we have:
$$\phi(e)=f_0(S_e)=f(S_e)-f(S_e)=0,$$
i.e. $\phi$ is supported in $X_b$, and minimizing $\|\phi\|_1$ is equivalent to minimizing $\|\phi-\phi_0\|_1$.\epr

\section{Algorithmic aspects in the case of trees}

In all of this section, $X=(V,E,\ell)$ denotes a metric tree with $|E|\ge 1$, and a function $f\in\mathcal{V}$ will be fixed.

\subsection{Computing $\|f\|_{TCE}$ with few additions}

Given $f\in \mathcal{F}(V)$, computing $\phi_f$ and then $\|f\|_{TCE}=\|\phi_f\|_1$\ involves many additions (of the order of $|E|\,|V|/2$). This can be improved.

\begin{Prop}\label{additions} Let $X=(V,E,\ell)$ be a metric tree with $|E|\geq 1$, and $f\in \mathcal{F}(V)$. 
One may compute $\|f\|_{TCE}=\|\phi_f\|_1$ with at most $2|E|-1$ additions.
\end{Prop}

\bpr We proceed by induction on $|E|$, the case $|E|=1$ being trivial. For $|E|>1$, let $v$ be a leaf (= terminal vertex) of $X$ and the corresponding edge $e$ is oriented with $e_+=v$. Then $\phi_f(e)=f(v)$. Delete the vertex $v$ and the edge $e$ and solve the problem for the remaining tree $Y$ and $f_1(w):=f(w)$ for $w\ne e_-$ and $f_1(e_-):=f(v)+f(e_-)$ (one addition). By induction hypothesis, $\|\phi_{f_1}\|_1$ can be computed in at most $2|E|-3$ additions, and $\|\phi_{f}\|_1=\|\phi_{f_1}\|_1+|f(v)|$.\epr

\subsection{Finding a transportation with minimal number of transports}

Let $f\in \mathcal{F}(V)$. We assume that we have already found $\phi_f$ (with $f=\partial (\phi_f)$). 

\begin{quest}\label{algo} Find a solution of STC involving the minimal number of transports (namely $\le |Supp(f)|-1$, see remark \ref{rem3.4}). 
\end{quest}

In other words, we want to find $P\in \tilde{\Pi}_f$ with $W(P)=\|f\|_{TC}$ and $|Supp(P)|\le |Supp(f)|-1$.

\begin{Rem}\label{rem4.3}
If $\phi_f(e)=0$ for some edges, we replace $E$ by the support $E'$ of $\phi_f$. Then $(V,E')$ is a forest. Treating separately each connected component, we will find an optimal transportation plan whose support has at most $|Supp(f)|-k$ elements, where $k$ is the number of connected components that contain points of the support of $f$.
\end{Rem}

Let $f\in \mathcal{F}(V)$. Using Remark \ref{rem4.3}, we may assume that $\phi_f(e)\ne 0$ for every $e\in E$. It is convenient for what follows to change the orientation of all edges with $\phi_f(e)<0$ and so assume that $\phi_f(e)> 0$ for every edge.

\medskip
Then a solution $P$ to Question \ref{algo} can be constructed algorithmically as follows.
\begin{enumerate}
\item Choose arbitrarily $a\in V_+$.
\item Construct a finite sequence of edges $e_1,\ldots,e_k$ such that $(e_1)_-=a$, $(e_j)_+=(e_{j+1})_-$ for $1\le j<k$, and $(e_k)_+\in V_-$.

To show the existence of this sequence, observe that $f(a)=\sum_{e,e_-=a}\phi_f(e)-\sum_{e,
e_+=a}\phi_f(e)$. It follows that there exists $e_1\in E$, with $(e_1)_-=a$. Assume that $e_1,\ldots,e_j$ are constructed,  since $0<f((e_j)_+)+\phi_f(e_j)\le \sum_{e,e_-=(e_j)_+}\phi_f(e)$ , we see that if $(e_j)_+\not\in V_-$, there exists $e_{j+1}\in E$ with $(e_j)_+=(e_{j+1})_-$.

\item Setting $b=(e_k)_+$ and $\alpha:=\min\{f(a),-f(b),\min\{\phi_f(e_j):1\le j\le k\}\}$, construct $P\in \tilde{\Pi}_f$ with $P(a,b)=\alpha$, $W(P)=\|f\|_{TC}$ and $|Supp(P)|\le |Supp(f)|-1$.



To do this, we proceed by induction on the number of elements of $Supp(f)$.

If $Supp(f)$ consists only of the two points $a,b$, then there is only one possible transportation plan with support equal to $V_+\times V_-$, given by
$$P(x,y)=\left\{\begin{array}{ccc}f(a) & \mbox{if} & (x,y)=(a,b) \\0 & \mbox{otherwise.} & \end{array}\right.$$

Assume then that $Supp(f)$ has $\ge3$ points. Define \begin{itemize}
\item $\psi:E\to \R_+$ by setting  $\psi(e)=\begin{cases}\alpha &$if $e \in \{e_1,\ldots,e_k\}\\0&$otherwise$\end{cases}$.
\item $h=-\partial \psi$ so that $h(a)=\alpha$, $h(b)=-\alpha$, $h(x)=0$ for $x\not\in \{a,b\}$.
\item $g=f-h$.
\end{itemize}

Define the plan $P_0$, by setting $P_0(a,b)=\alpha$ and $P_0(x,y)=0$ if $(x,y)\ne (a,b)$. Then, as in the initial step, $P_0\in \Pi_{h}$, and $\|\psi\|_1=W(P_0)$

Since $\psi$ and $\phi_f-\psi$ are nonnegative, $\|\phi_f\|_1=\|\psi\|_1+\|\phi_f-\psi\|_1$, in other words $\|f\|_{TC}=\|g\|_{TC}+\|h\|_{TC}$. We now consider 3 cases.
 
 \begin{enumerate}
\item If $\alpha =f(a)$, then $Supp(g)\subset Supp(f)\setminus \{a\}$. By the induction hypo-thesis, there is a transportation plan $Q\in \Pi_g$ with $W(Q)=\|g\|_{TC}$ and whose support has at most $|Supp(g)|-1$ elements. As $a\not\in Supp(g)$, $Q(a,b)=0$.

\item If $\alpha=-f(b)$ then $Supp(g)\subset Supp(f)\setminus \{b\}$: this case is similar to the previous one.

\item If $\alpha<\min\{f(a),-f(b)\}$, then $Supp(g)=Supp(f)$. There exists $j$ with $\alpha =\phi_f(e_j)$. Now $Supp(\phi_g)\subset E\setminus\{e_j\}$, and $(V,Supp(\phi_g))$ has at least two connected components, each of which has a strictly smaller number of vertices $v$ with $g(v)\ne 0$ (as $a$ and $b$ have been disconnected). By the induction hypothesis,  there is a transportation plan $Q\in \tilde{\Pi}_g$ with $W(Q)=\|g\|_{TC}$ and whose support has at most $|Supp(g)|-k$ elements - where $k$ is the number of connected components that contain points of $Supp(g)=Supp(f)$. Again, $Q(a,b)=0$.
\end{enumerate}

In all cases, $P=P_0+Q$ does the job, i.e. $P\in\tilde{\Pi}_f$, $P(a,b)=P_0(a,b)$,  
$$W(P)=W(P_0)+W(Q)=\|h\|_{TC}+\|g\|_{TC}=\|f\|_{TC},$$
and $|Supp(P)|\le|Supp(f)|-1$.\epr
\end{enumerate}

\begin{Rem}
Denote by $(V_f,E_P)$ the graph whose vertices are the points of the support of $f$ and edges are the points of the support of $P$. One may easily check that, in the above construction, if $(V_g,E_Q)$ is a forest, then $(V_f,E_P)$ is also a forest. In other words, our construction always yields a forest.\end{Rem}






\noindent
Authors addresses:\\

\noindent
Universit\'e Paris Cit\'e and Sorbonne Universit\'e\\ 
CNRS, IMJ-PRG\\
 F-75013 Paris, France\\
georges.skandalis@imj-prg.fr

\bigskip
\noindent
Institut de math\'ematiques\\
Universit\'e de Neuch\^atel\\
11 Rue Emile Argand - Unimail\\
CH-2000 Neuch\^atel - Switzerland\\
alain.valette@unine.ch

\end{document}